\input amstex
\documentstyle{amsppt}\nologo\footline={}

\def\PU{\mathop{\text{\rm PU}}}
\def\Area{\mathop{\text{\rm Area}}}
\def\S{\mathop{\text{\rm S}}}
\def\tr{\mathop{\text{\rm tr}}}
\def\SU{\mathop{\text{\rm SU}}}
\def\B{\mathop{\text{\rm B}}}
\def\G{\mathop{\text{\rm G}}}
\def\SL{\mathop{\text{\rm SL}}}
\def\Arg{\mathop{\text{\rm Arg}}}
\def\Re{\mathop{\text{\rm Re}}}

\hsize450pt

\topmatter\title A Hyperelliptic View on Teichm\"uller Space.
II\endtitle\author Sasha Anan$'$in\endauthor\thanks Partially supported
by the Institut des Hautes \'Etudes Scientifiques
(IH\'ES).\endthanks\address Departamento de Matem\'atica, IMECC,
Universidade Estadual de Campinas,\newline13083-970--Campinas--SP,
Brasil\endaddress\email Ananin$_-$Sasha\@yahoo.com\endemail\subjclass
30F60 (51M10, 57S30)\endsubjclass\abstract Using the methods of the
previous paper [ABG], we show that the Teichm\"uller space $\Cal T$ of
all closed Riemann surfaces is fibred twice over the Teichm\"uller
space $\Cal H$ of hyperelliptic ones. Both fibre bundles
$\pi_1,\pi_2:\Cal T\to\Cal H$ are real algebraic (rational). They
define an embedding $\Cal T\hookrightarrow\Cal H\times\Cal H$. In
addition, we indicate slight modifications of the proof of [ABG,
Theorem 5.1] providing an elementary proof of Toledo's rigidity
theorem.\endabstract\endtopmatter\document

\centerline{\bf1.~Introduction}

\medskip

It is always possible to {\it triangulate\/} a given group, i.e., to
choose generators $g_j$'s so that every defining relation takes the
form $g_1g_2g_3=1$. The well-known decomposition of a Riemann surface
into `pairs of pants' can be viewed as this kind of geometric
triangulation of the surface group. There are other geometric
triangulations:

Suppose that the surface group $G$ is a discrete subgroup in the group
$\Cal L$ of all orientation-preserving isometries of the hyperbolic
plane. Every $1\ne g\in G$ is known to be hyperbolic. Take the
isometries $1\ne g_j\in G$, $j=1,2$, and denote by $\G_j$ the axis of
$g_j$. If $\G_1$ and $\G_2$ intersect, $\G_1\cap\G_2=\{p_3\}$, then
$g_1=R(p_1)R(p_3)$ and $g_3=R(p_3)R(p_2)$ for suitable $p_1\in\G_1$ and
$p_2\in\G_2$, where $R(p)$ stands for the reflection in $p$. For
$g_3:=R(p_2)R(p_1)$, we have $g_1g_2g_3=1$. In such a manner, one can
easily create many triangulations of $G$. (In this respect, it seems
important to study all intersection points of axes of isometries in
$G$, i.e., all intersection points of closed geodesics in the
corresponding Riemann surface.) A hyperelliptic Riemann surface can be
easily characterized in these terms as a one that admits certain
generators $g_j$ (see the relations (2.1.1)) whose axes intersect all
in the same point. A suitable triangulation of $G$ gives raise to the
following global description of the Teichm\"uller space.

The Teichm\"uller space $\Cal T$ of closed Riemann surfaces of a fixed
genus $\ge2$ is fibred twice over the Teichm\"uller space $\Cal H$ of
hyperelliptic ones, $\pi_1,\pi_2:\Cal T\to\Cal H$. The bundles
$\pi_1,\pi_2$ are real algebraic (even rational, in the terms of the
Klein model) and induce an embedding
$\Cal T\hookrightarrow\Cal H\times\Cal H$. The fibres can be explicitly
described (see Theorem 2.3.5 and its proof).

At the end of this article, the reader may find a slight modification
of the proof of [ABG, Theorem~5.1] leading to an elementary and easy
proof of Toledo's rigidity theorem\footnote{Toledo's rigidity theorem
says that a representation $\varrho:G\to\PU(2,1)$ with maximal Toledo
invariant is faithful, discrete, and preserves a complex geodesic in
$\Bbb H_\Bbb C^2$.}
(Theorem 3.5). 

\medskip

{\bf Acknowledgements.} I am very grateful to Eduardo Carvalho Bento
Gon\c calves and Carlos Henrique Grossi Ferreira for their interest to
this work.

\bigskip

\centerline{\bf2.~Two fibre bundles over hyperelliptic Teichm\"uller
space}

\medskip

{\bf2.1.~Preliminaries.} As in [ABG], $W$ denotes a two-dimensional
$\Bbb C$-vector space equipped with a hermitian form of signature
$+-$ and $\Cal L:=\PU W$.

\smallskip

Let $n\ge6$ be an even integer. Denote by $H_n$ the group generated by
$r_1,r_2,\dots,r_n$ with the defining relations $r_i^2=1$,
$i=1,2,\dots,n$, and $r_n\dots r_2r_1=1$. The subgroup $G_n\le H_n$
constituted by the words of even length in the $r_i$'s is the surface
group $G_n=\pi_1\Sigma_g$, where $g=\frac n2-1$. As is easy to see,
$G_n$ admits the generators $g_i:=r_nr_i$, $i=1,2,\dots,n-1$, and the
defining relations
$$g_{n-1}g_{n-2}^{-1}g_{n-3}\dots g_2^{-1}g_1=1,\qquad g_{n-1}^{-1}
g_{n-2}g_{n-3}^{-1}\dots g_2g_1^{-1}=1.\leqno{\bold{(2.1.1)}}$$
Indeed, we can reconstruct $H_n$ from the group $G$ given by the
defining relations (2.1.1) as
$$H=\langle G,r\mid r^2=1,\ g_i^r=g_i^{-1},\ i=1,2,\dots,n-1\rangle$$
so that $r_n:=r$ and $r_i:=rg_i$, $i=1,2,\dots,n-1$.

As in [ABG], $\Cal R^+H_n$ and $\Cal R^+G_n$ denote the spaces of
representations in $\Cal L$ with maximal area $(n-4)\pi$ and
$2(n-4)\pi$. By Goldman's theorem [ABG, Theorems 3.15 and 5.1], such
representations are faithful and discrete. Denote by $\Cal H_n^+$ and
$\Cal T_n^+$  the corresponding Teichm\"uller spaces. We will use the
well-known fact that $\Cal R^+G_n$ and (hence) $\Cal T_n^+$ are
connected. Of course, one can find an elementary proof of this fact in
the style of [ABG, Section 5].

Let $\varrho\in\Cal R^+G_n$. The following remark is trivial.

\medskip

{\bf2.1.2.~Remark.} {\sl Let $1\ne g,g'\in G_n$. If the isometries
$\varrho g$ and $\varrho g'$ have a common fixed point $p\in\S W$, they
share a common axis.}

\medskip

{\bf Proof.} As is easy to see, $\tr[\varrho g,\varrho g']=2$ in the
terms of $\SU W$. Since $\varrho G_n$ contains no parabolic isometries,
we obtain $[\varrho g,\varrho g']=1$
$_\blacksquare$

\medskip

The connectedness of $\Cal R^+G_n$ and Remark 2.1.2 imply that, for
given $g,g'\in G_n$ such that $[g,g']\ne1$, the axes of the isometries
$\varrho g$ and $\varrho g'$ either are ultraparallel or intersect in
a point in $\B W$ independently of the choice of
$\varrho\in\Cal R^+G_n$.

We will frequently use the following well-known and trivial facts
concerning the group $\Cal I$ of all isometries of a given full
geodesic $\G$. The group $\Cal I$ reminds a dihedral group. It is
generated by the reflections in points in $\G$ and possesses the
$1$-parameter subgroup $\Cal V\le\Cal I$ of index $2$ formed by
hyperbolic isometries. Hence, $\Cal V\simeq(\Bbb R,+)$ and, for all
$t\in\Cal I\setminus\Cal V$ and $v\in\Cal V$, we have $t^2=1$,
$v^t=v^{-1}$, and there exists a unique $s\in\Cal I\setminus\Cal V$
such that $v=ts$.

\medskip

{\bf2.2.~Construction of fibre bundles.} Let $\varrho\in\Cal R^+G_n$.
Since, for $i=2,3,\dots,n-1$, the axes of $\varrho g_1$ and
$\varrho g_i$ intersect in the case of a hyperelliptic $\varrho$, they
intersect for an arbitrary $\varrho\in\Cal R^+G_n$. Denote by $t_i$ the
reflection in this intersection point and by $\G$, the axis of
$\varrho g_1$. Define
$$t_1:=t_{n-1},\qquad s_n:=t_{n-1}.\leqno{\bold{(2.2.1)}}$$
For suitable reflections $s_i$, we have
$$\varrho g_i=t_is_i,\qquad i=1,2,\dots,n-1,$$
where $s_1$ is the reflection in a point in $\G$. The relations (2.1.1)
imply the relations
$$t_{n-1}s_{n-1}s_{n-2}t_{n-2}t_{n-3}s_{n-3}\dots s_2t_2t_1s_1=1,\qquad
s_{n-1}t_{n-1}t_{n-2}s_{n-2}s_{n-3}t_{n-3}\dots t_2s_2s_1t_1=1$$
which can be rewritten as
$$s_ns_{n-1}s_{n-2}t_{n-2}t_{n-3}\dots s_3s_2t_2t_1s_1=1,\qquad
s_ns_{n-1}t_{n-1}t_{n-2}s_{n-2}s_{n-3}\dots t_3t_2s_2s_1=1.$$
Equivalently,
$$s_ns_{n-1}s_{n-2}h_{n-2}s_{n-3}s_{n-4}h_{n-4}\dots
h_4s_3s_2h_2s_1=1,\qquad s_ns_{n-1}h_{n-1}s_{n-2}s_{n-3}h_{n-3}\dots
h_3s_2s_1=1,$$
where $h_i:=t_it_{i-1}$, $i=2,3,\dots,n-1$, is an isometry with the
axis $\G$ or the identity. Define

\bigskip

$\bullet$\quad$u_n:=u_{n-1}:=u_{n-2}:=1$,\qquad$u_{n-3}:=
u_{n-4}:=h_{n-2}$,\qquad$u_{n-5}:=u_{n-6}:=h_{n-2}h_{n-4}$,\qquad
$\dots$,

\smallskip

\centerline{$u_3:=u_2:=h_{n-2}h_{n-4}\dots h_4$\qquad$u_1:=
h_{n-2}h_{n-4}\dots h_2$;}

\medskip

$\bullet$\quad$v_n:=v_{n-1}:=1$,\qquad$v_{n-2}:=v_{n-3}:=
h_{n-1}$,\qquad$v_{n-4}:=v_{n-5}:=h_{n-1}h_{n-3}$,\qquad$\dots$,

\smallskip

\centerline{$v_2:=v_1:=h_{n-1}h_{n-3}\dots h_3$;}

\medskip

$\bullet$\quad $a_1:=u_1s_1$,\qquad$c_1:=v_1s_1$;\qquad$a_i:=
s_i^{u_i}$,\qquad$c_i:=s_i^{v_i}$\qquad for\qquad $i=2,3,\dots,n$.

\bigskip

\noindent
(Note that $a_1$ and $c_1$ are reflections in points in $\G$.) We
obtain the relations
$$a_na_{n-1}a_{n-2}\dots a_3a_2a_1=1,\qquad
c_nc_{n-1}c_{n-2}c_{n-3}\dots c_2c_1=1$$
and, hence, two representations $\varrho_1,\varrho_2:H_n\to\Cal L$.

It is easy to observe that $s_1,s_2,\dots,s_n$ and
$t_1,t_2,\dots,t_{n-1}$ depend {\bf algebraically} on
$\varrho\in\Cal R^+G_n$ (not~involving radicals, when using the Klein
model). So do the $h_i$'s, $a_i$'s, and $c_i$'s. In particular,
the~maps $\pi_1:\varrho\mapsto\varrho_1$ and
$\pi_2:\varrho\mapsto\varrho_2$ are continuous, implying, in view of
[ABG, Lemma 3.2], that the maps
$\pi_\varepsilon:\varrho\mapsto\Area\varrho_\varepsilon$,
$\varepsilon=1,2$, are constant. If $\varrho$ is hyperelliptic, then
$t_i=\varrho r_n$ and $s_i=\varrho r_i$ for all $i=1,2,\dots,n-1$.
Hence, $h_i=1$ for all $i=1,2,\dots,n$ and
$\varrho_\varepsilon=\varrho$, $\varepsilon=1,2$. We conclude that the
areas of the $\varrho_\varepsilon$'s are maximal. Thus, we have
constructed two maps $\pi_1,\pi_2:\Cal R^+G_n\to\Cal R^+H_n$. For~the
reason of continuity, the isometries $\varrho g_1$, $\varrho_1g_1$, and
$\varrho_2g_1$ share the same repeller and attractor. Also,
the~isometry $\varrho_1r_n=\varrho_2r_n$ is the reflection in the
intersection of the axes of $\varrho g_1$ and $\varrho g_{n-1}$.

\smallskip

Let us show that $\varrho$ can be reconstructed from $\varrho_1$ and
$\varrho_2$. It suffices to reconstruct $t_i$ for all $i=1,2,\dots,n-1$
because the $t_i$'s provide the $u_i$'s which, in turn, allow to find
the $s_i$'s from the known $a_i$'s. To this aim, we write
$$c_i={a_i}^{f_i},\qquad i=1,2,3,\dots,n,$$
where the $f_i$'s are hyperbolic isometries with the axis $\G$ (or the
identities). Clearly, $f_i:=v_iu_i^{-1}$ for $i=2,\dots,n$ and $f_1$ is
determined by the equality $f_1u_1s_1f_1^{-1}=v_1s_1$.

For $1\le i\le n-3$, we have
$$u_i=t_{n-2}t_{n-3}\dots t_{i+1}\quad\text{if}\quad
i\equiv0\mod2,\qquad u_i=t_{n-2}t_{n-3}\dots t_i\quad\text{if}\quad
i\not\equiv0\mod2,\leqno{\bold{(2.2.2)}}$$
$$v_i=t_{n-1}t_{n-2}\dots t_i\quad\text{if}\quad i\equiv0\mod2,\qquad
v_i=t_{n-1}t_{n-2}\dots t_{i+1}\quad\text{if}\quad
i\not\equiv0\mod2.\leqno{\bold{(2.2.3)}}$$
In view of $t_1=t_{n-1}$, we have
$u_1v_1=t_{n-2}t_{n-3}\dots t_2t_1t_{n-1}t_{n-2}t_{n-3}\dots
t_2=(t_{n-2}t_{n-3}\dots t_2)^2=1$.
This implies $f_1=v_1$ because $v_1u_1s_1v_1^{-1}=v_1s_1$ is equivalent
to $u_1v_1=1$. Thus,
$$f_1=t_{n-1}t_{n-2}\dots t_3t_2,\qquad f_{n-2}=t_{n-1}t_{n-2},\qquad
f_{n-1}=1,\qquad f_n=1,$$
$$f_i=t_{n-1}t_{n-2}\dots t_{i+1}t_it_{i+1}\dots t_{n-2},\qquad
i=2,\dots n-3.\leqno{\bold{(2.2.4)}}$$
Note that
$$f_1^{-1}f_2f_3^{-1}\dots
f_{n-3}^{-1}f_{n-2}=1.\leqno{\bold{(2.2.5)}}$$
Indeed,
$$f_2f_3^{-1}\dots f_{2k-2}f_{2k-1}^{-1}f_{2k}f_{2k+1}^{-1}\dots
f_{n-3}^{-1}f_{n-2}=$$
$$=\big[(t_{n-1}t_{n-2}\dots t_3t_2)(t_3\dots
t_{n-2})\big]\big[(t_{n-2}\dots t_4t_3)(t_4\dots
t_{n-2}t_{n-1})\big]\dots$$
$$\dots\big[(t_{n-1}t_{n-2}\dots t_{2k-1}t_{2k-2})(t_{2k-1}\dots
t_{n-2})\big]\big[(t_{n-2}\dots t_{2k}t_{2k-1})(t_{2k}\dots
t_{n-2}t_{n-1})\big]$$
$$\big[(t_{n-1}t_{n-2}\dots t_{2k+1}t_{2k})(t_{2k+1}\dots
t_{n-2})\big]\big[(t_{n-2}\dots t_{2k+2}t_{2k+1})(t_{2k+2}\dots
t_{n-2}t_{n-1})\big]\dots$$
$$\dots\big[(t_{n-2}t_{n-3})(t_{n-2}t_{n-1})\big][t_{n-1}t_{n-2}]=
t_{n-1}t_{n-2}\dots t_3t_2=f_1.$$

Suppose that the $a_i$'s and $c_i$'s are known. Then the $f_i$'s are
known. From (2.2.1), we find $t_1=t_{n-1}=s_n=a_n$. From the known
$t_{n-1}$ and $f_{n-2}=t_{n-1}t_{n-2}$, we find $t_{n-2}$. If we have
already found $t_1,t_{n-1},t_{n-2},\dots,t_{i+1}$, then we can find
$t_i$ from the known
$f_i=t_{n-1}t_{n-2}\dots t_{i+1}t_it_{i+1}\dots t_{n-2}$.

\medskip

{\bf2.3.~Description of fibres.} We go back along the above way of
constructing the maps $\pi_\varepsilon:\Cal R^+G_n\to\Cal R^+H_n$,
$\varepsilon=1,2$, and describe the fibres of these maps. We deal only
with $\pi_1$ (similar considerations work for $\pi_2$).

\smallskip

Let us fix a representation $\varrho_1\in\Cal R^+H_n$. For suitable
points $q_j\in\B W$, $j=1,2,\dots,n$, we have
$a_j:=\varrho_1 r_j=R(q_j)$ and
$$a_na_{n-1}\dots a_2a_1=1.\leqno{\bold{(2.3.1)}}$$
Denote by $\G$ the full geodesic $\G{\prec}q_n,q_1{\succ}$ and by
$\Cal I\ge\Cal V$ the group of all isometries of $\G$ and its subgroup
formed by the hyperbolic ones. We plan to describe the fibre
$\pi_1^{-1}\varrho_1\simeq\pi_2\pi_1^{-1}\varrho_1$ as formed by the
representations $\varrho_2\in\Cal R^+H_n$ given by a relation of the
form
$$a_na_{n-1}a_{n-2}^{f_{n-2}}a_{n-3}^{f_{n-3}}\dots
a_2^{f_2}a_1^{f_1}=1,\leqno{\bold{(2.3.2)}}$$
where $f_j\in\Cal V$, $i=1,2,\dots,n-3,n-2$, satisfy (2.2.5).

First, we do not take (2.2.5) into account. In other words, we look for
points $q'_1,q'_2,\dots,q'_{n-3},q'_{n-2}$ placed in the same curves
equidistant from $\G$ as the corresponding points
$q_1,q_2,\dots,q_{n-3},q_{n-2}$ such that the relation
$R(q_n)R(q_{n-1})R(q'_{n-2})R(q'_{n-3})\dots R(q'_2)R(q'_1)=1$ is valid
and provides a representation $\varrho_2$ with maximal area. Of course,
$q'_j=f_jq_j$ for suitable (and unique) $f_j\in\Cal V$,
$j=1,2,\dots,n-3,n-2$.

To this aim, we construct subsequently points $q'_j\in D_j$ and
simultaneously the positive $1$-cycle of $\varrho_2$ (see [ABG,
Definition 3.12]), where $q_j\in D_j$ stands for the curve equidistant
from $\G$, $j=n-2,n-3,\dots,2,1$. By [ABG, Lemma 3.11],
$q_{n-1},q_{n-2},q_{n-3},\dots,q_2$ are all on the side of the normal
vector to $\G$. So are the $D_j$'s. Initially, we have
$b^1,e^1\in\S W\cap\G$, the repeller and the attractor of~$a_1a_n$.
Hence, we obtain $b^{n-2}:=a_{n-1}b^1$ and $e^{n-2}:=a_{n-1}e^1$ such
that the cycle $b^1,e^1,b^{n-2},e^{n-2}$ is positive. By induction on
$4\le k\le n-2$, suppose that we have already constructed points
$q'_j\in D_j$ for $j=n-2,n-3,\dots, k+1$ and the positive cycle
$b^1,e^1,b^k,e^k,b^{k+1},e^{k+1},\dots,b^{n-3},e^{n-3},b^{n-2},
e^{n-2}$.
It is easy to see that, in order to have the positive cycle
$b^1,e^1,b^{k-1},e^{k-1},b^k,e^k,b^{k+1},e^{k+1},\dots,b^{n-3},e^{n-3},
b^{n-2},e^{n-2}$,
where $b^{k-1}:=R(q'_k)b^k$ and $e^{k-1}:=R(q'_k)e^k$, it is necessary
and sufficient to place $q'_k\in\B W$ in the open region $R$ on the
side of the normal vector to the geodesic $\G(b^k,e^1)$. Since
$\G(b^k,e^1)$ and $D_k$ have the common vertex $e^1$ and the other
vertex of $D_k$ does not belong to the closure of $R$, the~intersection
$\G(b^k,e^1)\cap D_k\subset\B W$ consists of a single point $d_k$. Now
we simply (have to) choose $q'_k$ in the open segment of $D_k$ between
the points $d_k$ and $e^1$. It is essential to observe that $d_k$
depends algebraically on $\varrho_1$ and on
$q'_{n-2},q'_{n-3},\dots,q'_{k+1}$. Finally, we construct
$q'_{n-2},q'_{n-3},\dots,q'_4$ and the positive cycle
$b^1,e^1,b^3,e^3,b^4,e^4,\dots,b^{n-3},e^{n-3},b^{n-2},e^{n-2}$.

\medskip

{\bf2.3.3.~Lemma.} {\sl Let\/ $b^1,e^1,b^3,e^3\in\S W$ be a positive
cycle and let\/ $D_2,D_3$ be curves equidistant from the geodesic\/
$\G(b^1,e^1)$ and situated on the side of the normal vector to\/
$\G(b^1,e^1)$. Then there exist unique\/ $p_j\in D_j$, $j=2,3$, such
that the cycle\/ $b^1,e^1,b^2,e^2,b^3,e^3$ is positive, where\/
$R(p_2)b^1=b^2$, $R(p_2)e^1=e^2$, $R(p_3)b^2=b^3$, and
$R(p_3)e^2=e^3$.}

\medskip

Using Lemma 2.3.3, we construct the points $q'_3\in D_3$ and
$q'_2\in D_2$. Then, we find a unique point $q'_1\in D_1=\G$ providing
the relation
$R(q_n)R(q_{n-1})R(q'_{n-2})R(q'_{n-3})\dots R(q'_2)R(q'_1)=1$ as in
the proof of [ABG, Corollary 3.17]. By [ABG, Theorem 3.15], we obtain a
representation $\varrho_2:H_n\to\Cal L$ with maximal area. In addition,
we have $q'_j=f_jq_j$ for suitable $f_j\in\Cal V$,
$j=1,2,\dots,n-3,n-2$.

\smallskip

Suppose that the $f_j$'s satisfy (2.2.5). Put $t_{n-1}:=t_1:=a_n$. From
$f_{n-2}=t_{n-1}t_{n-2}$, determine $t_{n-2}$. Using (2.2.4), define
subsequently the reflections $t_i$ for $i=n-3,n-2,\dots,3,2$. All the
$t_i$'s are reflections in points in $\G$. For $1\le i\le n-3$, define
$u_i$ and $v_i$ by means of (2.2.2) and (2.2.3). Put
$s_1:=u_1^{-1}a_1$, $s_i:=u_i^{-1}a_iu_i$ for $i=2,3,\dots,n-3$,
$s_{n-2}:=a_{n-2}$, $s_{n-1}:=a_{n-1}$. Finally, define
$$g'_i:=t_is_i,\qquad i=1,2,\dots,n-1.$$

The equalities (2.2.4) and the calculus after (2.2.5) show that
$f_1=t_{n-1}t_{n-2}\dots t_3t_2$. It follows from (2.2.2),
$t_1=t_{n-1}$, and  (2.2.3) that
$u_1f_1=t_{n-2}\dots t_3t_2t_1t_{n-1}t_{n-2}\dots t_3t_2=(t_{n-2}\dots
t_3t_2)^2=1$
and $f_1=v_2$. Thus, $v_2^{-1}f_1^2u_1=1$.

The relation $a_na_{n-1}\dots a_2a_1=1$ can be written as
$$t_{n-1}s_{n-1}s_{n-2}u_{n-3}s_{n-3}u_{n-3}^{-1}u_{n-4}s_{n-4}
u_{n-4}^{-1}u_{n-5}\dots
u_4^{-1}u_3s_3u_3^{-1}u_2s_2u_2^{-1}u_1s_1=1.$$
It easily follows from (2.2.2) that, for $1\le i\le n-4$,
$$u_{i+1}^{-1}u_i=1\quad\text{if}\quad i\equiv0\mod2,\qquad
u_{i+1}^{-1}u_i=t_{i+1}t_i\quad\text{if}\quad i\not\equiv0\mod2.$$
Hence, we obtain
$$t_{n-1}s_{n-1}s_{n-2}t_{n-2}t_{n-3}s_{n-3}s_{n-4}t_{n-4}t_{n-5}\dots
t_4t_3s_3s_2t_2t_1s_1=1,$$
that is,
$$g'_{n-1}{g'}_{n-2}^{-1}g'_{n-3}{g'}_{n-4}^{-1}\dots
g'_3{g'}_2^{-1}g'_1=1.$$

It follows from (2.2.2), (2.2.3), and (2.2.4) that $f_i=v_iu_i^{-1}$
for $i=2,3,\dots,n-3$. Therefore, in view of $f_{n-2}=t_{n-1}t_{n-2}$,
the relation
$a_na_{n-1}a_{n-2}^{f_{n-2}}a_{n-3}^{f_{n-3}}\dots
a_2^{f_2}a_1^{f_1}=1$
can be written as
$$t_1s_{n-1}t_{n-1}t_{n-2}s_{n-2}t_{n-2}t_{n-1}s_{n-3}^{v_{n-3}}
s_{n-4}^{v_{n-4}}\dots s_2^{v_2}f_1^2u_1s_1=1.$$
It easily follows from (2.2.3) that, for $1\le i\le n-4$,
$$v_{i+1}^{-1}v_i=t_{i+1}t_i\quad\text{if}\quad i\equiv0\mod2,\qquad
v_{i+1}^{-1}v_i=1\quad\text{if}\quad i\not\equiv0\mod2.$$
Since $v_{n-3}=t_{n-1}t_{n-2}$ and $v_2^{-1}f_1^2u_1=1$, we obtain
$$s_{n-1}t_{n-1}t_{n-2}s_{n-2}s_{n-3}t_{n-3}t_{n-4}\dots
t_3t_2s_2s_1t_1=1,$$
that is,
$${g'}_{n-1}^{-1}g'_{n-2}{g'}_{n-3}^{-1}\dots g'_2{g'}_1^{-1}=1.$$
Thus, we constructed a representation $\varrho:G_n\to\Cal L$ given by
$\varrho: g_i\mapsto g'_i$, $i=1,2,\dots,n-1$.

\medskip

{\bf2.3.4.~Lemma.} {\sl Let\/ $\varrho\in\Cal R^+H_6$. Put\/
$a_i:=\varrho r_i$, $1\le i\le6$, and denote by\/ $\Cal V$ the\/
$1$-parameter group of hyperbolic isometries that contains\/ $a_6a_1$.
The representations in\/ $\Cal R^+H_6$ corresponding to the relations
of the type\/ $a_6a_5a_4^{f_4}a_3^{f_3}a_2^{f_2}a_1^{f_1}=1$, where\/
$f_1,f_2,f_3,f_4\in\Cal V$, form a space that is isomorphic to\/
$\Cal V$ by means of the map\/
$(f_1,f_2,f_3,f_4)\mapsto f_1^{-1}f_2f_3^{-1}f_4$.}

\medskip

Let $\varrho_1\in\Cal R^+H_n$ be given by (2.3.1). We are going to
show that the $(n-2)$-tuples
$(f_1,f_2,\dots,\allowmathbreak f_{n-3},f_{n-2})$ in (2.3.2) that
provide $\varrho_2\in\Cal R^+H_n$ and satisfy (2.2.5) form a connected
space $F_{\varrho_1}$.

When, in the considerations before Lemma 2.3.3, we choose points
$q'_{n-2},q'_{n-3},\dots,q'_6,q'_5$ such that the cycle
$b^1,e^1,b^4,e^4,$ $b^5,e^5,\dots,b^{n-3},e^{n-3},b^{n-2},e^{n-2}$ is
positive, where $R(q_{n-1})b^1=b^{n-2}$, $R(q_{n-1})e^1=e^{n-2}$,
$R(q'_j)b^j=b^{j-1}$, and $R(q'_j)e^j=e^{j-1}$ for
$j=n-2,n-3,\dots,6,5$, we choose in fact a point
$(f_5,f_6,\dots,f_{n-3},f_{n-2})$ given by the relations $q'_j=f_jq_j$,
$j=5,6,\dots,n-3,n-2$, in a space diffeomorphic to $\Bbb R^{n-6}$.

Denote $p_6:=q_n$ and $\{p_5\}:=\G(b^1,b^4)\cap\G(e^1,e^4)$. Let
$q_j\in D_j$, $j=2,3,4$, stand for the curves equidistant from the
geodesic $\G:=\G(b^1,e^1)$ and situated on the side of the normal
vector to $\G$. Choose a point $p_4:=f'_4q_4\in D_4$ on the side of the
normal vector to $\G(b^4,e^1)$. Define $R(p_4)b^4=b^3$ and
$R(p_4)e^4=\nomathbreak e^3$. So, the cycle $b^1,e^1,b^3,e^3,b^4,e^4$
is positive. By Lemma 2.3.3, find $p_2:=f'_2q_2\in D_2$,
$p_3:=f'_3q_3\in D_3$, and, finally, $p_1:=f'_1q_1\in\G$ such that the
relation $R(p_6)R(p_5)R(p_4)R(p_3)R(p_2)R(p_1)=1$ is valid and provides
a representation $\varrho\in\Cal R^+H_6$. By Lemma 2.3.4, there exist
unique $f''_1,f''_2,f''_3,f''_4$ such that
$p_6,p_5,f''_4p_4,f''_3p_3,f''_2p_2,f''_1p_1$ provide a
representation in $\Cal R^+H_6$ and
$${f''}_1^{-1}f''_2{f''}_3^{-1}f''_4=
f'_1{f'}_2^{-1}f'_3{f'}_4^{-1}f_5f_6^{-1}\dots f_{n-3}f_{n-2}^{-1}.$$
So, we can find unique $f_j:=f''_jf'_j$, $j=1,2,3,4$, such that (2.2.5)
holds and (2.3.2) provides a representation $\varrho_2\in\Cal R^+H_n$.
In other words, $F_{\varrho_1}\simeq\Bbb R^{n-6}$.

Since $F_{\varrho_1}$ is connected and $\Area\varrho$ depends
continuously on $\varrho$, we conclude that the representation
$\varrho$ constructed above from the representations $\varrho_1$ and
$\varrho_2$ that are given by (2.3.1) and (2.3.2) belongs
to~$\Cal R^+G_n$. (It suffices to take $f_j=1$ for all
$j=1,2,\dots,n-3,n-2$.)

\medskip

{\bf2.3.5.~Theorem.} {\sl There exist two algebraic fibre bundles\/
$\pi_1,\pi_2:\Cal T_n^+\to\Cal H_n^+$ that define an embedding\/
$\Cal T_n^+\hookrightarrow\Cal H_n^+\times\Cal H_n^+$. The fibre\/
$\pi_1^{-1}[\varrho]\simeq\pi_2\pi_1^{-1}[\varrho]\simeq\Bbb R^{n-6}$
is the space of all\/ $(n-2)$-tuples\/ $(f_1,f_2,\dots,f_{n-2})$ of
hyperbolic isometries in\/ $\PU W$ that share the axis of\/ $a_na_1$,
meet\/ {\rm(2.2.5)}, and satisfy the relation\/
$a_na_{n-1}a_{n-2}^{f_{n-2}}a_{n-3}^{f_{n-3}}\dots
a_2^{f_2}a_1^{f_1}=1$
providing a representation with maximal area, where\/
$a_j:=\varrho r_j$.\/ {\rm(}A~similar fact is valid for\/
$\pi_2$.{\rm)} $_\blacksquare$}

\medskip

{\bf2.4.~Proofs of Lemmas 2.3.3 and 2.3.4.} Working in the upper
half-plane model, we write
$\infty:=\left[\smallmatrix 1\\0\endsmallmatrix\right]$ and
$p:=\left[\smallmatrix p\\1\endsmallmatrix\right]$ for $p\in\Bbb C$. A
curve $D$ equidistant from the geodesic $\G:=\G(\infty,0)$ and situated
on the side of the normal vector to $\G$ is simply a ray
$D=(1+ik^{-1})\Bbb R^+$ with $k>0$. The reflection in the point
$p=(1+ik^{-1})t\in D$ is given by the matrix
$$R(p)=\left[\matrix-k&(k+k^{-1})t\\-kt^{-1}&k\endmatrix\right]\in
{\SL}_2\Bbb R,\qquad R(p)x=t-\frac{t^2}{k^2(x-t)},\qquad
R(p)t=\infty.$$

For $q>0$, the intersection $D\cap\G(q,0)$ is the point
$\displaystyle\frac{(1+ik^{-1})q}{1+k^{-2}}$. Hence, the part
of $D$ that lies on the side of the normal vector to $\G(q,0)$
can be parameterized as
$$p(u):=(1+ik^{-1})t(u),\qquad t(u):=\frac q{(1+k^{-2})(1+u)},\qquad
u>0.\leqno{\bold{(2.4.1)}}$$
By a straightforward calculus, we obtain
$$R\big(p(u)\big)x=\frac{q\big((1+u)x-q\big)}
{(1+u)\big((1+k^{-2})(1+u)x-q\big)}.\leqno{\bold{(2.4.2)}}$$

\medskip

{\bf Proof of Lemma 2.3.3.} Let $D_j=(1+ik_j^{-1})\Bbb R^+$ for
$j=2,3$, $b^1:=\infty$, and $e^1:=0<b^3<e^3$. By~(2.4.1), every point
$p_3\in D_3$ on the side of the normal vector to $\G(b^3,e^1)$ has the
form $p_3=(1+ik_3^{-1})t_3$, where
$t_3=\displaystyle\frac{b^3}{(1+k_3^{-2})(1+u)}$ and $u>0$. Put
$b^2:=R(p_3)b^3$ and $e^2:=R(p_3)e^3$. Then $0<b^2<e^2<b^3$.
By (2.4.2),
$$b^2=\frac{b^3u}{(1+u)(k_3^{-2}+u+k_3^{-2}u)},\qquad
e^2=\frac{b^3\big(e^3(1+u)-b^3\big)}
{(1+u)\big(e^3(1+k_3^{-2})(1+u)-b^3\big)}.$$
The fact that $R\big((1+ik_2^{-1})t_2\big)b^2=b^1$ means that
$t_2=b_2$. Therefore, we have $p_j=(1+ik_j^{-1})t_j$, $j=2,3$, where
$$t_2=\frac{b^3u}{(1+u)(k_3^{-2}+u+k_3^{-2}u)},\qquad
t_3=\frac{b^3}{(1+k_3^{-2})(1+u)}.\leqno{\bold{(2.4.3)}}$$
The fact that $R\big((1+ik_2^{-1})t_2\big)e^2=e^1$ means that
$e^2=(1+k_2^{-2})t_2$. So, we obtain the following condition for $u>0$
$$\frac{b^3\big(e^3(1+u)-b^3\big)}{(1+u)\big(e^3(1+k_3^{-2})(1+u)-
b^3\big)}=(1+k_2^{-2})\frac{b^3u}{(1+u)(k_3^{-2}+u+k_3^{-2}u)}.$$
This condition is equivalent to the equation
$$(1+k_3^2)u^2+u+\Big(1-\frac{b^3}{e^3}\Big)(k_3^2u-k_2^2u-k_2^2)=0
\leqno{\bold{(2.4.4)}}$$
that possesses a unique solution $u>0$ because $e^3>b^3>0$
$_\blacksquare$

\medskip

{\bf Proof of Lemma 2.3.4.} We assume that the repeller and attractor
of $a_1a_6$ are $b^1=\infty$ and $e^1=0$. Put $b:=a_5b^1$ and
$e:=a_5e^1$. The cycle $b^1,e^1,b,e$ is positive. This means that
$0<b<e$. It is immediate that
$a_5=\displaystyle\frac1{\sqrt{be-b^2}}\left[\smallmatrix-b&be\\-1&b
\endsmallmatrix\right]$.
For suitable $q_j\in\B W$ and $k_j>0$, we have $a_j=R(q_j)$ and
$q_j\in d_j\Bbb R^+$, where $d_j:=1+ik_j^{-1}$, $j=2,3,4$. Also,
$a_6=R(d_1)$ for some $d_1\in\G(b^1,e^1)$.

The space of the representation in Lemma 2.3.4 can be described as
follows. Take an arbitrary $p_4\in d_4\Bbb R^+$ on the side of the
normal vector to $\G(b,e^1)$. Define $b^3:=R(p_4)b$ and $e^3:=R(p_4)e$.
Clearly, the~cycle $b^1,e^1,b^3,e^3,b,e$ is positive. By Lemma 2.3.3,
find $p_j\in d_j\Bbb R^+$, $j=2,3$. There exists $p_1\in\G(b^1,e^1)$
such that the relation $a_6a_5R(p_4)R(p_3)R(p_2)R(p_1)=1$ is valid and
provides a representation in $\Cal R^+H_6$. Note that, for every choice
of $p_4$, the corresponding points $p_3,p_2,p_1$ are unique.

Every isometry $f\in\Cal V$ has the form
$f=\left[\smallmatrix r&0\\0&r^{-1}\endsmallmatrix\right]$, $r>0$. Such
an $f$ acts as a dilatation: $fx=r^2x$. For suitable $t_j>0$, we have
$p_j=d_jt_j$, $j=1,2,3,4$. Therefore, it suffices to show that the map
$(p_1,p_2,p_3,p_4)\mapsto t_1t_2^{-1}t_3t_4^{-1}$ is an isomorphism.

By (2.4.1) and (2.4.2), for an arbitrarily chosen $v>0$, we have
$p_4=d_4t_4$,
$$t_4=\frac b{(1+k_4^{-2})(1+v)},\
b^3=\frac{bv}{(1+v)(k_4^{-2}+v+k_4^{-2}v)},\
e^3=\frac{b\big(e(1+v)-b\big)}
{(1+v)\big(e(1+k_4^{-2})(1+v)-b\big)},\leqno{\bold{(2.4.5)}}$$
$t_2$ and $t_3$ are given by (2.4.3) where $u>0$ is a unique solution
of (2.4.4).

The reflection in the point $id$, $d>0$, is given by the matrix
$R(id)=\left[\smallmatrix0&d\\-d^{-1}&0\endsmallmatrix\right]
\in\SL_2\Bbb R$.
Therefore,
$a_6R(p_1)=R(d_1)R(d_1t_1)=\left[\smallmatrix-t_1^{-1}&0\\0&-t_1
\endsmallmatrix\right]$
and
$t_1=-\left[\smallmatrix0&1\endsmallmatrix\right]a_6R(p_1)
\left[\smallmatrix0\\1\endsmallmatrix\right]$.
In the terms of $\SL_2\Bbb R$, we~have
$a_6R(p_1)=\pm a_5R(p_4)R(p_3)R(p_2)$. Using (2.4.3) and (2.4.5), we
obtain
$$w:=\pm\frac{k_2(k_3+k_3^{-1})\sqrt{be-b^2}}
{k_4+k_4^{-1}}t_1t_2^{-1}t_3t_4^{-1}=$$
$$-\frac{k_2(k_3+k_3^{-1})}{k_4+k_4^{-1}}t_2^{-1}t_3t_4^{-1}
\left[\smallmatrix0&1\endsmallmatrix\right]\left[\smallmatrix-b&be\\
-1&b\endsmallmatrix\right]\left[\smallmatrix-k_4&(k_4+k_4^{-1})t_4\\
-k_4t_4^{-1}&k_4\endsmallmatrix\right]\left[\smallmatrix-k_3&
(k_3+k_3^{-1})t_3\\-k_3t_3^{-1}&k_3\endsmallmatrix\right]
\left[\smallmatrix-k_2&(k_2+k_2^{-1})t_2\\-k_2t_2^{-1}&k_2
\endsmallmatrix\right]\left[\smallmatrix0\\1\endsmallmatrix\right]$$
$$=\frac{(1+u+k_3^2u)v}{1+u}.$$
Then $u=\displaystyle\frac{w-v}{(1+k_3^2)v-w}$. It follows from (2.4.5)
that
$1-\displaystyle\frac{b^3}{e^3}=
\frac{(e-b)(v+1)}{(k_4^2v+v+1)(ev+e-b)}$.
In the terms of $v$ and $w$, the equation (2.4.4) takes the form
$$\frac{e(1+k_4^2)}{e-b}vw(w-v)+\Big(\frac e{e-b}+k_4^2\Big)w(w-v)+
(1+k_2^2+k_3^2)(v+1)(w-v)-k_2^2k_3^2v(v+1)=0.\leqno{\bold{(2.4.6)}}$$
Since $e>b>0$ and $k_2,k_3\ne0$, the equation (2.4.6) admits a unique
solution $w(v)>v$ for every $v>0$. (Such solution corresponds to the
solution $u>0$ of (2.4.4).) In particular, $w(v)\to+\infty$ as
$v\to+\infty$. As~is easy to see, $w(v)\to+0$ as $v\to+0$. Consider the
real plane cubic $C$ given by (2.4.6). It intersects the infinite
projective line at three points corresponding to the vertical,
horizontal, and diagonal straight lines. If the smooth function
$v\mapsto w(v)$ is not a diffeomorphism of $(0,+\infty)$, it has
a local maximum $w(v_0)$ followed by a local minimum $w(v_1)$,
$0<v_0<v_1$. Hence, the horizontal line $H$ passing through the point
$p:=\big(v_1,w(v_1)\big)$ is tangent to $C$ at $p$. On the other hand,
$H$ intersects $C$ in some finite point $\big(v,w(v)\big)$ with
$0<v<v_0$. This contradicts B\'ezout's Theorem
$_\blacksquare$

\bigskip

\centerline{\bf3.~Elementary proof of Toledo's rigidity theorem}

\medskip

We assume here that the reader is familiar with [ABG, Section 2], with
[ABG, Definition 3.6], with the second part of [ABG, Remark 3.7], with
[ABG, Section 5] from the beginning till Lemma 5.13, including the
proofs, and with a certain part of the proof of [AGG, Proposition
2.1.6]. We follow the line given in [ABG]. The difference is that $W$
is now a three-dimensional $\Bbb C$-vector space equipped with a
hermitian form of signature $++-$.

\smallskip

Let $u,p_1,p_2,p_3\in\B W$. As is shown in the proof of [AGG,
Proposition 2.1.6], the integral of the K\"ahler potential $P_u$ along
the geodesic $\G[p_1,p_2]$ equals
$\displaystyle\int_{\G[p_1,p_2]}P_u=-\frac\pi2+\frac12\Arg\frac{\langle
u,p_2\rangle\langle p_2,p_1\rangle}{\langle u,p_1\rangle}$.
It~follows that the integral of the K\"ahler form over the oriented
geodesic triangle equals
$$\int_{\Delta(p_1,p_2,p_3)}\omega=\int_{\partial\Delta(p_1,p_2,p_3)}
P_u=-\frac{3\pi}2+\frac12\sum\limits_{i=1}^3\Arg\frac{\langle
u,p_{i+1}\rangle\langle p_{i+1},p_i\rangle}{\langle u,p_i\rangle}$$
(the indices are modulo $3$). Taking $u=p_1$, we get
$\displaystyle\int_{\Delta(p_1,p_2,p_3)}\omega=
-\frac\pi2+\frac12\Arg(g_{13}g_{32}g_{21})$,
where $g_{ij}:=\langle p_i,p_j\rangle$. Now the formula
$$\Area\Delta(p_1,p_2,p_3):=-\int_{\Delta(p_1,p_2,p_3)}\omega=\frac12
\arg(-g_{12}g_{23}g_{31})\leqno{\bold{(3.1)}}$$
easily follows. Indeed, all we need to show is that
$\Re(-g_{12}g_{23}g_{31})\ge0$. By Sylvester's criterion,
$\det[g_{ij}]\le0$, which is equivalent to
$$2\Re\frac{g_{12}g_{23}g_{31}}{g_{11}g_{22}g_{33}}\ge
\frac{g_{12}g_{21}}{g_{11}g_{22}}+\frac{g_{23}g_{32}}{g_{22}g_{33}}+
\frac{g_{31}g_{13}}{g_{33}g_{11}}-1$$
because $g_{11},g_{22},g_{33}<0$. It remains to observe that
$\displaystyle\frac{g_{12}g_{21}}{g_{11}g_{22}}\ge1$.

The formula (3.1) is similar to [ABG, (2.1)] and introduces the `area'
of the oriented geodesic triangle $\Delta(p_1,p_2,p_3)$,
$p_1,p_2,p_3\in\overline\B W$, having no coinciding isotropic vertices.
Define also $\Area\Delta(p,p,q):=0$ for $p\in\S W$ and
$q\in\overline\B W$. From now on, we follow [ABG, Section 2] till the
following assertion:

\medskip

{\bf3.2.~Remark.} {\sl$\Area(c;p_1,p_2,p_3)=\Area\Delta(p_1,p_2,p_3)$
for all pairwise distinct\/ $c,p_1,p_2,p_3\in\S W$.}

\medskip

{\bf Proof.} Varying the pairwise distinct points
$c,p_1,p_2,p_3\in\S W$, we can reach the situation where they all
belong to a complex geodesic
$_\blacksquare$

\medskip

{\bf3.3.~Definition.} We say that the points
$p_1,p_2,\dots,p_k\in\S V$, $k\ge3$, form a {\it positive cycle\/} if
they all belong the same complex geodesic and form there a positive
cycle in the sense of [ABG, Definition 3.6]. Note that by this
definition the points $p_1,p_2,\dots,p_k$ are pairwise distinct.

\medskip

The second part of [ABG, Remark 3.7] implies immediately the

\medskip

{\bf3.4.~Remark.} {\sl If the cycles\/ $p_1,p_2,\dots,p_k\in\S W$,
$k\ge3$, and\/ $p_k,p_{k+1},p_1\in\S W$ are positive, then the cycle\/
$p_1,p_2,\dots,p_k,p_{k+1}$ is positive.}

\medskip

Now we pass to [ABG, Section 5], define the group $G_n$ for even
$n\ge6$, define $\Cal L:=\PU W$, and state~the

\medskip

{\bf3.5.~Theorem {\rm[Tol]}.} {\sl A representation\/
$\varrho:G_n\to\Cal L$ is faithful and discrete if\/
$\Area\varrho=\pm\displaystyle\frac{(n-4)\pi}2$. In this case, there
exists a stable complex geodesic\/ $C$, $\varrho G_nC=C$.}

\medskip

{\bf Proof.} We proceed literally by [ABG, Section 5] till [ABG, Remark
5.10]. The introduced area $\Area\varrho$ of a representation $\varrho$
is exactly (one fourth of) the Toledo invariant. This can be easily
seen from the considerations at the beginning of this section and from
those at the very beginning of [Tol]. Now~we should substitute the
words `hyperbolic isometry' by `loxodromic isometry' and continue
proceeding by [ABG, Section 5]. After we have proven [ABG, Lemma 5.13],
we are done. Indeed, since the cycle $s_3,t_3,s_4,t_4$ is positive, by
shifting the indices, we conclude that all the $s_i$'s and $t_i$'s
belong to the same complex geodesic. It remains to observe that $G_n$
is generated by the $h_i$'s and apply Goldman's theorem [ABG, Theorem
5.1]
$_\blacksquare$

\bigskip

\centerline{\bf4.~References}

\medskip

[ABG]~S.~Anan$'$in, E.~C.~B.~Gon\c calves, {\it A Hyperelliptic View
on Teichm\"uller Space. {\rm I},} available at
http://arxiv.org/abs/0709.1711

\smallskip

[AGG]~S.~Anan$'$in, C.~H.~Grossi, N.~Gusevskii, {\it Complex Hyperbolic
Structures on Disc Bundles over Surfaces. {\rm I.}~General Settings. A
Series of Examples,} available at http://arxiv.org/abs/math/0511741

\smallskip

[Tol]~D.~Toledo, {\it Representations of surface groups in complex
hyperbolic space,} J.~Differential Geom. {\bf29} (1989), No.~1,
125--133.

\enddocument